\documentclass[twoside,leqno,twocolumn]{article}

\usepackage[letterpaper]{geometry}
\usepackage{ltexpprt}
\usepackage{hyperref}
\usepackage{cite}
\usepackage{graphicx}
\usepackage{amsfonts}
\usepackage{amsmath,bm}
\usepackage{enumitem}
\usepackage[mathscr]{euscript}
\usepackage{cancel}
\usepackage{caption}
\usepackage{subcaption}
\usepackage{tikz}
\usepackage{pgfplots}
\usetikzlibrary{spy,calc}
\usepackage{hyperref}
\usepackage{setspace}
\usepackage{float}
\usepackage{titlesec}
\usepackage{wrapfig}
\usepackage{fancyhdr}
\begin{document}

\title{\Large Control-Affine Extremum Seeking Control with Attenuating Oscillations: A Lie Bracket Estimation Approach \thanks{Dr. Eisa acknowledges the fund of the 2023 Office of Research URC Faculty Scholars Research Award at the University of Cincinnati, OH-USA.}}
\author{Sameer Pokhrel\thanks{Ph.D. student at the Department of Aerospace Engineering and Engineering Mechanics, University of Cincinnati, OH, USA. Email: pokhresr@mail.uc.edu}
\and Sameh A. Eisa \thanks{Assistant professor at the Department of Aerospace Engineering and Engineering Mechanics, University of Cincinnati, OH, USA. Email: eisash@ucmail.uc.edu and sameheisa235@hotmail.com.}}

\date{}

\maketitle


\begin{abstract} \small \baselineskip=9pt
Control-affine Extremum Seeking Control (ESC) systems have been increasingly studied and applied in the last decade. In a recent effort, many control-affine ESC structures have been generalized in a unifying class and their stability was analyzed. However, guaranteeing vanishing oscillations at the extremum point for said class requires strong conditions that may not be feasible or easy to check/design by the user, especially when the gradient of the objective function is unknown. In this paper, we introduce a control-affine ESC structure that remedies this problem such that: (i) its oscillations attenuate structurally via a novel application of a geometric-based Kalman filter and a Lie bracket estimation approach; and (ii) its stability is characterized by a time-dependent (one-bound) condition that is easier to check and relaxed when compared to the generalized approach mentioned earlier. We provide numerical simulations of three problems to demonstrate the effectiveness of our proposed ESC; these problems cannot be solved with vanishing oscillations using the generalized approach in the literature.
\end{abstract}
\section{Introduction} \label{s: Intro}
Extremum Seeking Control (ESC) is a model-free adaptive control technique that stabilizes a steady state map of a dynamical system around the extremum point of an objective function that we have access to its measurements, but not its expression. In 2000, what is known as the classic ESC structure \cite{KRSTICMain} has been analyzed through singular perturbation and averaging (refer \cite{Maggia2020higherOrderAvg} for more on averaging theory), and its stability was characterized. This classic structure has been extended to many other forms and found many applications \cite{ariyur2003real,scheinker2017modelfreebook}. 

In this paper, we focus on control-affine ESC systems \cite{ESCTracking, DURR2013, scheinker2016boundedDither,BoundedUpdateKrstic,ExpStabSUTTNER2017,eisa2023}, which often deal with problems/applications naturally expressed in control-affine formulation -- such as but not limited to multi-agent systems. In control-affine ESC structures, the objective function is not perturbed/modulated directly by the input signal, rather it -- or a functional of it -- is multiplied by the input signal, resulting in the objective function being incorporated within the vector fields of the control-affine ESC system. This has invoked the use of Lie bracket-based analysis as these tools are very natural to control-affine systems as known in geometric control theory \cite{bullo2019geometric}. Durr et al. \cite{DURR2013} first used the concept of Lie Bracket System (LBS) approximation of ESCs, mainly for stability characterization. They called it ``the corresponding LBS." They showed that under certain assumptions, the practical stability property of a control-affine ESC system follows from the asymptotic stability property of the corresponding LBS; in essence, LBS characterization of control-affine ESCs is analogous to singular perturbation and averaging characterization of classic ESC-related structures. In the following years to \cite{DURR2013}, many researchers have adopted structures that utilize LBS-based approaches for characterizing stability and improving the ESC performance (e.g.,  \cite{DURR2017,GRUSHKOVSKAYA2017}). Recently, Grushkovskaya et al. \cite{VectorFieldGRUSHKOVSKAYA2018} provided a generalized class for many control-affine ESC structures \cite{scheinker2016boundedDither,BoundedUpdateKrstic,ExpStabSUTTNER2017} including the work of Durr et al. \cite{DURR2013}. In their work, certain conditions and bounds if guaranteed on the objective function and its derivatives, one should be able to guarantee vanishing control inputs at the extremum point (this effort has been extended in \cite{taha2021vanishing}), hence asymptotic convergence of the ESC system is achieved. 

Since Durr et al. \cite{DURR2013} introduced their ESC approach, until the generalization of control-affine ESCs by Grushkovskaya et al. \cite{VectorFieldGRUSHKOVSKAYA2018}, there seems to be a consistent trade-off: approaches of ESC with simple stability condition but persistent oscillations (e.g., \cite{DURR2013}), or approaches of ESC with complex stability conditions to guarantee vanishing oscillations (e.g., the generalized ESC approach in \cite{VectorFieldGRUSHKOVSKAYA2018}). Our motivation is to balance said trade-off by achieving relaxed stability conditions with vanishing oscillations. By examining the literature of classic ESC-related structures, the contribution made recently in \cite{AttenuatedOscillaiion2021} has made it possible to implement the classic ESC structure with attenuating oscillations and without over-complicating, the stability condition of the classic ESC structure in \cite{KRSTICMain}. The key was finding the proper adaptation law for the amplitude of the control input. We aim at introducing an analogous concept to  \cite{AttenuatedOscillaiion2021} \textit{but for control-affine ESC systems} to achieve the results of \cite{VectorFieldGRUSHKOVSKAYA2018} with minimal complications to the relaxed stability conditions in \cite{DURR2013}. In this paper, we propose a LBS-estimation approach for control-affine ESC systems which: (i) have attenuating oscillations that vanish at the extremum point; and (ii) its stability is characterized by one time-dependent condition/bound, with no bounds on the gradient derivatives. 
Moreover, the framework proposed in this paper introduces a novel adaptation law utilizing a novel geometric-based Kalman filtering (GEKF) approach to estimate the LBS approximating a control-affine ESC. We solve numerically three problems including a multi-agent system to demonstrate the effectiveness of the proposed approach; examples chosen in this paper are for problems the generalized ESC approach in \cite{VectorFieldGRUSHKOVSKAYA2018} cannot solve with vanishing oscillations.
%
%
%
%
%
\section{Background and preliminaries}
\begin{figure}[ht]
    \centering
\includegraphics[width=0.45\textwidth]{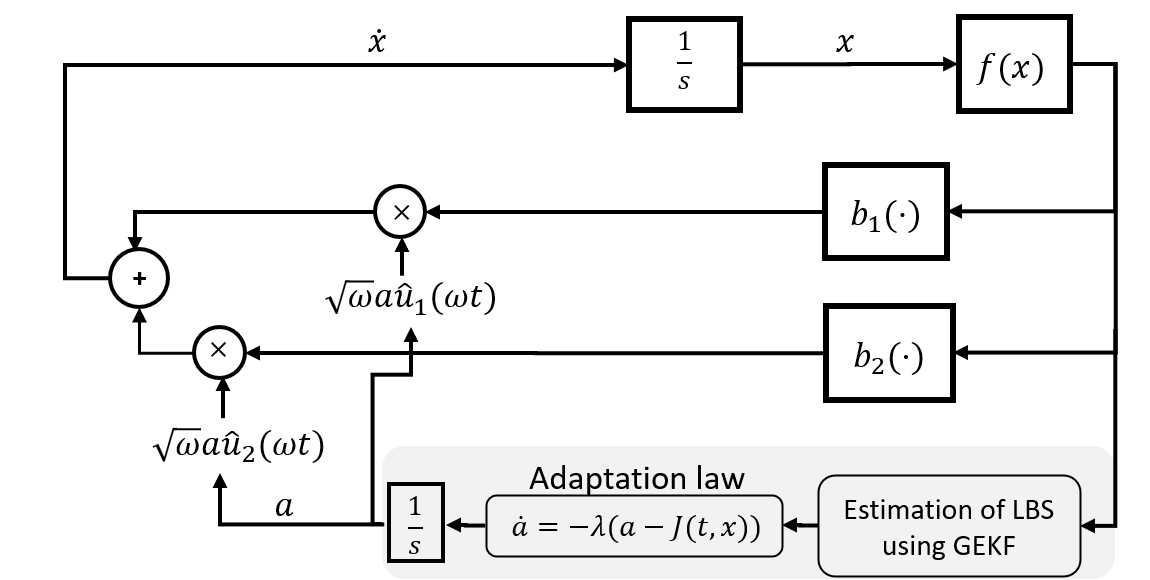}
    \caption{Proposed ESC structure.}
    \label{fig:ESC_scheme}
\end{figure}

Control-affine ESC systems can be written as \cite{DURR2013}:
\begin{equation}\label{eqn:ESC_back}
    \dot{\bm{x}}=\bm{b_d}(t,\bm{x})+ \sum\limits_{i=1}^{m} \bm{b_i}(t,\bm{x})\sqrt{\omega}u_i(t,\omega t),\\
\end{equation}
with $\bm{x}(t_0)=\bm{x_0}\in \mathbb{R}^n$ and $\omega \in (0,\infty)$.
Here $\bm{x}$ is the state space vector, $\bm{b}_d$ is the drift vector field of the system, $u_i$ are the control inputs, $m$ is the number of control inputs, and $\bm{b}_i$ are the control vector fields.
Now, the corresponding LBS \cite{DURR2013} to (\ref{eqn:ESC_back}) is defined as in \eqref{eqn:Lie_back}:
\begin{equation}\label{eqn:Lie_back}
    \dot{\bm{z}}=\bm{b_d}(t,\bm{z})+ \sum_{\substack{i=1\\j=i+1}}^m [\bm{b_i},\bm{\bm{b_j}}](t,\bm{z})\nu_{j,i}(t),
\end{equation}
with $\nu_{j,i}(t)=\frac{1}{T}\int_0^T u_j(t,\theta)\int_0^\theta u_i(t,\tau)d\tau d\theta$. 
The operation $[\cdot,\cdot]$ is the Lie bracket operation between two vector fields $\bm{b_i},\bm{b_j}: \mathbb{R} \times\mathbb{R}^n \rightarrow \mathbb{R}^n$ with $\bm{b_i}(t,\cdot),\bm{b_j}(t,\cdot)$ being continuously differentiable, and is defined as $[\bm{b_i},\bm{b_j}](t,\bm{x}):=\frac{\partial \bm{b_j}(t,\bm{x})}{\partial \bm{x}}\bm{b_i}(t,x)-\frac{\partial \bm{b_i}(t,\bm{x})}{\partial \bm{x}}\bm{b_j}(t,\bm{x})$. A particular class of the control-affine ESC in (\ref{eqn:ESC_back}), which generalizes many ESC systems in control affine form found in literature \cite{ESCTracking, DURR2013, scheinker2016boundedDither,BoundedUpdateKrstic,ExpStabSUTTNER2017,VectorFieldGRUSHKOVSKAYA2018, taha2021vanishing}, is provided below:
\begin{equation}\label{eqn:ESC}
    \dot{{x}}={b_1}(f({x}))u_1+ {b_2}(f({x}))u_2, 
\end{equation}
where $f({x})$ is the objective function. This generalized structure of ESCs \eqref{eqn:ESC} is shown in the unshaded parts of figure \ref{fig:ESC_scheme}, with $u_1=a\sqrt{\omega}\hat{u}_1(\omega t),u_2=a \sqrt{\omega}\hat{u}_2(\omega t) $, $a \in \mathbb{R}$ is the amplitude of the input signal. Now, we impose the following assumptions on ${b}_1,{b}_2$, $\hat{u}_1$ and $\hat{u}_2$.
\begin{enumerate}[label=A\arabic*.]
    \item 
    $b_i\in C^2: \mathbb{R} \to \mathbb{R}, i={1,2}$ and for a compact set $\mathscr{C} \subseteq \mathbb{R}$, there exist $A_1, ..., A_3 \in [0,\infty)$ such that $|b_i(x)|\leq A_1,
    | \frac{\partial b_i(x)}{\partial x}|\leq A_2,
    |\frac{\partial [{b_j},{b_k}](x)}{\partial x}| \leq A_3$ for all $x\in \mathscr{C}, i={1,2};\; j={1,2};\; k={1,2}.$
    \item 
    $\hat{u}_i: \mathbb{R} \times \mathbb{R} \to \mathbb{R} , i=1,2$, are measurable functions. Moreover, there exist constants $M_i \in (0,\infty) $ that $sup_{\omega t \in \mathbb{R}}|\hat{u}_i(\omega t)|\leq M_i$, and 
    $\hat{u}_i(\cdot)$ is T-periodic, i.e. $\hat{u}_i(\omega t + T)=\hat{u}_i(\omega t),$ and has zero average, i.e. $\int_0^T \hat{u}_i(\tau) d\tau = 0,$ with $T \in (0,\infty)$ for all $\omega t \in \mathbb{R}$.
    \item 
    There exists an ${x}^* \in \mathscr{C}$ such that $\nabla f({x}^*)=0, \nabla f({x})\ne 0$ for all ${x}\in \mathscr{C}\backslash \{{x}^*\}; f({x}^*)=f^* \in \mathbb{R}$ is an isolated extremum value. 
\end{enumerate}

\begin{Remark}\label{remark:assumptions1_3}
Assumptions A1-A3 are typical in ESC literature \cite{KRSTICMain,DURR2013}. A1-A2 ensure well-posedness, boundedness, zero-mean average, and measurability of the control inputs. A3 ensures the objective function $f$ has an isolated local extremum $f^*$ at $\bm{x}^*$. 
\end{Remark}


\begin{theorem}\label{thm:from2018}
    (from \cite{VectorFieldGRUSHKOVSKAYA2018}) Let A1-A2 be satisfied, then for the ESC system in (\ref{eqn:ESC}), the corresponding LBS is: 
    \begin{equation}\label{eqn:Lie}
    \dot{{z}}=-\nu_{2,1} \nabla f({z}) b_0(f({z})),
    \end{equation}
    where $b_0(z)=b_2(z) \frac{db_1(z)}{dz}-b_1(z)\frac{db_2(z)}{dz}, z\in \mathbb{R}$.
    For the multi-variable ESC ($\bm{x}\in \mathbb{R}^n$) in the following form
\begin{equation}\label{eqn:ESC_multi}
    \dot{\bm{x}}=\sum \limits_{i=1}^n \left( b_{1i}(f({\bm{x}}))u_{1i}+ {b_{2i}}(f({\bm{x}}))u_{2i}\right) e_i, 
\end{equation}
the corresponding LBS is given by
\begin{equation}\label{eqn:Lie_multi}
    \dot{\bm{z}}=-\sum \limits_{i=1}^n \nu_{2i,1i} \frac{\partial f(\bm{z})}{\partial z_i} b_{0i}(f({\bm{z}})) e_i.
\end{equation}
where $f:  \mathbb{R}^n \to \mathbb{R}$, $e_i$ denotes the $i^{th}$ unit vector in $\mathbb{R}^n$.
\end{theorem}

Moreover, we introduce Theorem \ref{thm:esc_lbs} from \cite{DURR2013} (this theorem is also stated in \cite{VectorFieldGRUSHKOVSKAYA2018}) which links the stability properties of LBSs (\ref{eqn:Lie}) and (\ref{eqn:Lie_multi}) with ESC systems (\ref{eqn:ESC}) and (\ref{eqn:ESC_multi}), respectively. 

\begin{theorem}\label{thm:esc_lbs}
Let assumptions A1-A3 be satisfied and suppose that a compact set $\mathscr{C}$ is locally (uniformly) asymptotically stable for (\ref{eqn:Lie_multi}). Then $\mathscr{C}$ is locally practically (uniformly) asymptotically stable for (\ref{eqn:ESC_multi}).
\end{theorem}

Finally, we introduce the Chen-Fliess functional expansion \cite[chapter 3]{isidori1985nonlinear} that can be used to describe representations of input-output behavior of a nonlinear system of the control-affine form (\ref{eqn:ESC_back}) and associated output function $y=h(\bm{x})$, $y\in \mathbb{R}$. Let $T$ be a fixed value of the time and $u_1,...,u_m$ are real-valued piecewise continuous functions defined on [0,$T$]. Now, iterated integral for each multi-index ($i_k,...,i_0$) is defined as: 
\begin{equation}
    \int_0^t d\xi_{i_k}...d\xi_{i_0}=\int_0^t d\xi_{i_k}(\tau) \int_0^t d\xi_{i_{k-1}}... \int_0^t d\xi_{i_0},
\end{equation}
where $0\le t \le T$; $\xi_0(t)=t$; $\xi_i(t)=\int_0^t u_i(\tau)d\tau$ for $1 \le i \le m$. Now, the evolution of the output $y(t)$ for a short time $t\in[0,T]$ is given by: 
\begin{equation}\label{eqn:chenFliess}
    y(t)=h(\bm{x}_0)+\sum \limits_{k=0}^\infty \sum \limits_{i_0 , ..., i_k = 0}^m  L_{\bm{b}_{i_0}}... L_{\bm{b}_{i_k}} h(\bm{x}_0) \int_0^t d\xi_{i_k}...d\xi_{i_0},
\end{equation}
where $L_{\bm{b}}h$ is the Lie derivative of $h$ along $\bm{b}$.

\section{Main Results and Discussion} \label{s: MainResults}
The LBS in (\ref{eqn:Lie_multi}) approximates, captures the behavior of, and characterizes the ESC system in (\ref{eqn:ESC_multi}). One can easily observe that, with the exception of the gradient components, the LBS in  (\ref{eqn:Lie_multi}) is obtainable by access to measurements of the objective function $f(\bm{x})$ and the predetermined structure choices. For instance, $b_0$ per (\ref{eqn:Lie}) is determined by $b_1$ and $b_2$ built in the structure, and the term $\nu_{2,1}$ is dependent on the chosen control inputs. As a result, the LBS estimation is dependent on the gradient of the objective function. This leads us to formally introduce the concept of estimated LBS, denoted as $\hat{\bm{z}}$, and defined as follows:
\begin{equation}\label{eqn:LieEst}
    \dot{\hat{\bm{z}}}=-\sum \limits_{i=1}^n ( \nu_{2i,1i} \frac{\partial f(\bm{z})}{\partial z_i} b_{0i}(f({\bm{z}})) + \eta_i(t)) e_i =\bm{J}(t,\bm{z}).
\end{equation}
The estimated LBS in (\ref{eqn:LieEst}) is similar to (\ref{eqn:Lie_multi}) but with the addition of an error term ${\eta}_i(t)$ representing the inaccuracies introduced during the estimation. we impose the following assumption on ${\eta}_i(t)$:
\begin{enumerate}[label=A\arabic*.]
    \setcounter{enumi}{3}
    \item
     ${\eta}_i(t): \mathbb{R} \rightarrow \mathbb{R},i=1,...,n$ is measurable function and there exist constants $\theta_0,\epsilon_0 \in (0,\infty)$ such that $|{\eta}_i(t_2)-{\eta}_i(t_1)| \le \theta_0|t_2-t_1|$ for all $t_1,t_2 \in \mathbb{R}$ and $\sup_{t \in \mathbb{R}}|{\eta}_i(t)| \le \epsilon_0$.
    Furthermore, $\lim_{t\to \infty} {\eta}_i (t)={0}$.
\end{enumerate}
Assumption A4 implies that the error term $\eta_i(t)$ is measurable, bounded, and decreases with time. Having the knowledge of the estimated LBS, we now propose an ESC structure that couples the control-affine ESC in (\ref{eqn:ESC_multi}) with an adaptation law for the amplitude of the control input which depends on the estimated LBS as:
\begin{align}
    \dot{\bm{x}}&=\sum \limits_{i=1}^n \left( b_{1i}(f({x}))\sqrt{\omega} a_i(t)\hat{u}_{1i}+ {b_{2i}}(f({x}))\sqrt{\omega} a_i(t) \hat{u}_{2i}\right) e_i \label{eqn:generalizedSystem},\\
    \dot{\bm{a}}&=\sum \limits_{i=1}^n \left(-\lambda_i ({a}_i(t)-{J}_i(t,\bm{x})) \right ) e_i,\label{eqn:law}
\end{align}
where $a_i \in \mathbb{R}$ is the amplitude of the input signal, and $\lambda_i >0\in \mathbb{R}$ is a tuning parameter. Our proposed structure for $n=1$ is shown in figure \ref{fig:ESC_scheme}. Note that our structure in figure \ref{fig:ESC_scheme} without the shaded area is the same as the generalized ESC systems in \cite{VectorFieldGRUSHKOVSKAYA2018}. 
%
Next, we introduce our theorem on the proposed structure.

\begin{theorem}\label{thm:proposed_theorem}
Let A1-A4 be satisfied with some $\omega \in (\omega^*,\infty), \omega^* >0$ and suppose $\exists \; t^* >t_0 =0$ such that $\forall t>t^*$  and $|J_i(t,\bm{z})| \le 1/t^p$ with some $p>1$ then (i) the equilibrium point $\hat{\bm{z}}^* \in \mathscr{C}$ is locally asymptotically stable for the estimated LBS in \eqref{eqn:LieEst}, (ii) $a_i$ in (\ref{eqn:law}) is asymptotically convergent to 0; and (iii) the system in (\ref{eqn:generalizedSystem}) is practically asymptotically stable.\\
\emph{The proof can be found in the appendix \ref{sec:Proof}}.
\end{theorem}

\begin{Remark} Theorem \ref{thm:proposed_theorem} mean that if after some time the estimated right-hand side of the LBS is consistently bounded via $1/t^p$ ($p>1$), then the proposed ESC system is stable with guaranteed vanishing oscillations.
\end{Remark}

The above-mentioned results, until this point, depend on the knowledge of the estimated LBS. So, it is necessary to introduce a sound technique to obtain the estimated LBS inline with theorem 3.1. Estimating the gradient of the objective function is not a new concept in the ESC literature; however, it was done only to classic ESC-related structures using Kalman Filtering (KF) \cite{Chicka2006,AttenuatedOscillaiion2021,ESC_KF}. In said structures, KF measurement update equation was easily found via natural Taylor expansion of the perturbed objective function \cite{Chicka2006,AttenuatedOscillaiion2021}. That is, the addition of the input $du$ after the integrator in such structures perturbs the function argument as: 
\begin{equation}\label{eqn:taylor}
    f(x)=f(\hat{x}+du)=f(\hat{x})+f'(\hat{x}) du+O(|du|^2),
\end{equation}
where $\hat{x}$ denotes current estimate of the variable $x$ and $f'(\hat{x})$ denotes the gradient estimate.
However, (\ref{eqn:taylor}) \textit{cannot be used for control affine ESC structures} (the focus of this work) as the input signal $du$ does not perturb the function argument in such systems. It is important to emphasize that there is hardly  \textit{any} literature for estimating the gradient or LBSs in control-affine ESCs. Here, we introduce a novel and applicable geometric-based KF framework to control-affine ESC systems which we call "GEKF." The Chen-Fliess series (\ref{eqn:chenFliess}) (see also \cite{grushkovskaya2021extremum}) can be used to analyze control-affine ESC systems. Arguably, the Chen-Fliess series is a control-affine analog to the Taylor series. We set $y=f(x)$ in (\ref{eqn:chenFliess}) and truncate after first-order terms:
\begin{multline}\label{eqn:measurementEqn1_mainresult}
    f(\bm{x})|_{t_2}= f(\bm{x})|_{t_1}+ {\bm{b}_1} \cdot \nabla f(\bm{x})|_{t_1} U_1 +{\bm{b}_2} \cdot \nabla f(\bm{x}) |_{t_1} U_2 \\
    + O((\Delta t)^2),
\end{multline}
with $U_1 = \int_{t_1}^{t_2}{u}_1 d\tau, U_2 = \int_{t_1}^{t_2}{u}_2 d\tau$, $ t_2=t_1+\Delta t$, $\Delta t$ is the time period/step. We set $\Delta t$ to be a very small step size when compared to the periodic time. Hence, $\Delta t = K/\omega$, with some constant $K$, and the higher order terms are in $O((1/\omega)^2)$. The higher-order terms are sufficiently small when the frequency chosen is sufficiently large. 
Now, (\ref{eqn:measurementEqn1_mainresult}) can be used as the measurement update equation for our introduced continuous-discrete extended KF (GEKF) to estimate the LBS. It is worth noting that in our structure (figure \ref{fig:ESC_scheme}), the measurement of the objective function is taken from the ESC system. So, a zero mean error is expected to persist in the measurement equation given the oscillatory nature of the ESC trajectory about the averaged LBS. We assume the sum of this error and the higher order term $O(1/\omega)^2)$ in (\ref{eqn:measurementEqn1_mainresult}) as Gaussian measurement noise $\nu (t) \sim N(0,r)$. Next, we formulate the state propagation model of the GEKF to approximate the LBS corresponding to a control-affine ESC system with state space $x \in \mathbb{R}^n$. We define the state variables of the GEKF as follows:
\begin{align}\label{eqn:GE_states}
    \bm{\bar{X}}=\begin{bmatrix}
    [\bar{\bm{x}}_1]_{n \times 1}\\
    [\bar{\bm{x}}_2]_{n \times 1}\\
    [\bar{x}_3]_{1 \times 1}\\
    \end{bmatrix}=\begin{bmatrix}
    -\sum \limits_{i=1}^n \nu_{2i,1i} \frac{\partial f(\bm{x})}{\partial x_i} b_{0i}(f({\bm{x}})) e_i \\
    \dot{\bar{\bm{x}}}_1\\
     f(\bm{x})|_{t_1}\\
    \end{bmatrix},
\end{align}
where $\bar{\bm{x}}_1$ is equivalent to the right-hand side of the LBS (\ref{eqn:Lie_multi}). We assume $\bar{\bm{x}}_1$ has a constant derivative given by $\bar{\bm{x}}_2$. The state $\bar{x}_3$ is taken as the constant value of the objective function. Thus, the state propagation model, used in between the measurements at $t_1$ and $t_2$ is:
\begin{align}\label{eqn:GEKFDynamics}
    \dot{\bar{\bm{X}}}=\begin{bmatrix}
    \bar{\bm{x}}_2,
    0,
    0\\
    \end{bmatrix}^T + \bm{\Omega},
\end{align}
where $\bm{\Omega}$ is a random variable representing the process noise which is assumed to be Gaussian with zero mean and covariance $\bm{Q}$. 
 Furthermore, the process noise and the measurement noise are assumed to be uncorrelated. Note that the first $n$ states of the filter are estimating the right-hand side of the LBS in the filter dynamics. Hence the filter state itself is the integration of that which is proportional to the average. So, by the time the ESC reaches the limit cycle, the average of said right-hand side will incorporate the exact gradient. This can be extracted from the states of the filter to finally obtain the estimated LBS.
 The use of the GEKF to estimate the right-hand side of the LBS, as discussed above, introduces an error. We argue that, due to the choice of sufficiently large frequency and choice of our propagation dynamics and measurement equation of the filter, the error is bounded and decreases with time, i.e. it has the same characteristics as the error $\eta(t)$ introduced in (\ref{eqn:LieEst}) and satisfies assumption A4. Thus, our proposed technique of estimation can be used to obtain the estimated LBS.

It is essential to briefly highlight the advantages of the proposed control-affine ESC structure when compared to the generalized approach of Grushkovskaya et al. \cite{VectorFieldGRUSHKOVSKAYA2018} which unified most significant control-affine ESC systems  \cite{ESCTracking, DURR2013, scheinker2016boundedDither,BoundedUpdateKrstic,ExpStabSUTTNER2017}. First highlight: our proposed ESC stability is characterized by a time-dependent condition that, in principle, is simpler and easier to verify for a given ESC system when compared to the stability conditions introduced in \cite{VectorFieldGRUSHKOVSKAYA2018}. The stability conditions of \cite{VectorFieldGRUSHKOVSKAYA2018} -- provided in the appendix \ref{sec:StabilityConditionsB} -- are harder to obtain through estimation/approximation methods applied to the measurements of the objective function as said conditions require bounds on the objective function, its gradient, and even higher order derivatives.
Second highlight: in order to guarantee vanishing oscillations of the control inputs at the extremum point, a strong condition from \cite{VectorFieldGRUSHKOVSKAYA2018} (provided in appendix \ref{sec:StabilityConditionsB} as B2) has to be satisfied. However, B2 is hard to verify, check, or apply, especially without knowing the objective function or its derivatives as discussed above. Moreover, there is no apparent alternative for how one can attenuate the oscillations if B2 is not satisfied. Our approach on the other hand guarantees the vanishing of oscillations even when condition B2 is not satisfied. Next, we show the merit of our approach by studying three ESC cases. The first case is a simple single-variable ESC system taken from \cite{VectorFieldGRUSHKOVSKAYA2018}, the second case is a multi-variable single-agent vehicle case taken from \cite{BoundedUpdateKrstic}, and the third case is the multi-agent problem using single-integrator dynamics similar to Durr et al. \cite{DURR2013}. In the first and second cases, oscillations can not be attenuated as stated in \cite{VectorFieldGRUSHKOVSKAYA2018} due to condition B2 not being satisfied. We also show in appendix \ref{sec:StabilityConditionsB} that the third case problem from \cite{DURR2013} cannot be solved with vanishing oscillation due to the violation of B2. Nevertheless, with our proposed approach, the three ESC cases are solved with vanishing oscillations where the generalized approach \cite{VectorFieldGRUSHKOVSKAYA2018} did not succeed.

\section{Simulation Results}\label{s:Simulation}
\begin{figure}[ht]
    \centering
    \includegraphics[width=0.4\textwidth]{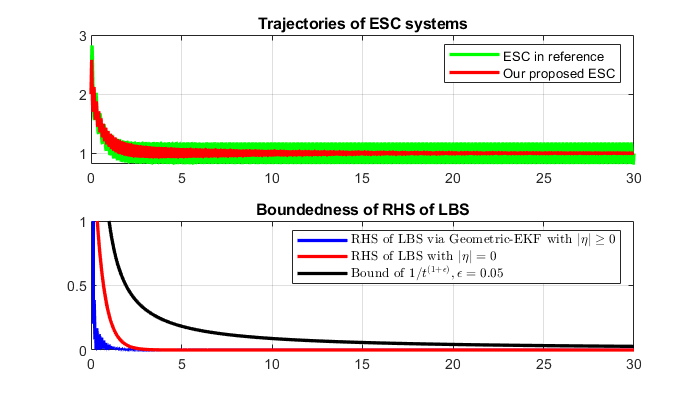}
    \caption{Comparison between our ESC system (convergent) with that of  \cite{VectorFieldGRUSHKOVSKAYA2018} (persistently oscillating) (top), simulation verification of Theorem \ref{thm:proposed_theorem} (bottom).}
    \label{fig:comparison_case1}
\end{figure}

\begin{figure}[ht]
    \includegraphics[width=0.5\textwidth, keepaspectratio]{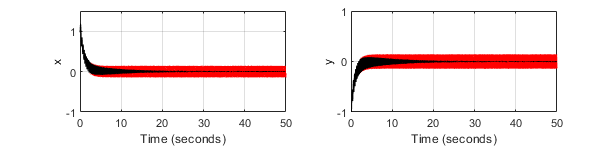}
    \caption{Trajectories of x and y coordinates employing ESC from \cite{BoundedUpdateKrstic} (red) vs. our proposed system (black).}
    \label{fig:plotWithTime_BU}
\end{figure}

\begin{figure}[ht]
    \centering
    \includegraphics[width=0.46\textwidth, keepaspectratio]{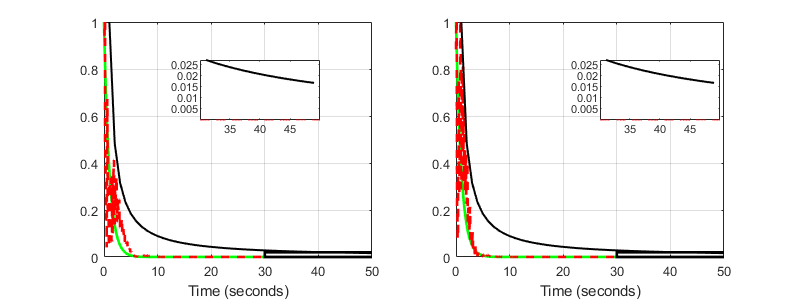}
    \caption{Verification of Theorem \ref{thm:proposed_theorem}: $|J_{x}|$ (red, left), $|J_{x}|$ with $|\eta| = 0 $ (green, left), $|J_{y}|$ (red, right), $|J_{y}|$ with $|\eta|=0$ (green, right), and $1/t^{1+0.05}$ (black) vs. time.}
    \label{fig:RHS_LBS_BU}
\end{figure}

\textbf{Case 1.}
This simple ESC system was used in \cite{VectorFieldGRUSHKOVSKAYA2018} to demonstrate that with the failure to satisfy condition B2 (refer the appendix), the ESC will have non-vanishing oscillations at the extremum point. We resolved the same problem by our proposed ESC (\ref{eqn:generalizedSystem})- (\ref{eqn:law}), and it works effectively. The equation of the system is $\dot{x}=f(x)u_1(t)+ u_2 (t)$ with $f(x)=2(x-x^*)^2, x\in \mathbb{R}, x^*=1$, $u_1=a\sqrt{\omega}\cos(\omega t)$ $u_2=a\sqrt{\omega}\sin(\omega t), \omega=8$. The LBS is $\dot{z}=- \alpha/2 \nabla f(z)$ with $\alpha = a^2$, and the expression for estimated LBS is $\dot{z}=- \alpha/2 \nabla f(z)+\eta (t)$. We used the GEKF as described in section \ref{s: MainResults}. Figure \ref{fig:comparison_case1} shows the simulation results for this case with the initial condition as $a_0=1, x_0=2$ and parameter $\lambda =0.1$ with a top plot showing the advantage of our proposed ESC (vanishing oscillation). Similarly, verification of Theorem \ref{thm:proposed_theorem} is provided in the bottom plot of figure \ref{fig:comparison_case1} where the estimated and exact right-hand side of the LBS shown are under the viable bound $1/t^p,p=1.05$.
\\
\noindent\textbf{Case 2.} In this case, we consider a single vehicle example moving in a two-dimensional plane in GPS-denied environment from \cite{BoundedUpdateKrstic}. The ESC system drives the vehicle to its optimal position based on measurable, but analytically unknown, cost function $f(x,y)$, where $x$ and $y$ are the coordinates of the vehicle. This problem as mentioned in \cite{VectorFieldGRUSHKOVSKAYA2018} cannot be solved with vanishing oscillations as condition B2 is not satisfied.  The equation of motion for the system is given as
\begin{equation}\label{eqn:ESC_BU}
\begin{split}
    \dot{x}= \cos({k f}) u_1 -  \sin({k f}) u_2\\
    \dot{y}=  \sin({k f}) u_1 + \cos({k f}) u_2,
\end{split}
\end{equation}
where $\omega$ is the amplitude and angular frequency of the signal and $k$ is a constant. The control inputs are taken as $u_1 = a \sqrt{\omega}\cos \omega t$ and $u_2=a\sqrt{\omega} \sin \omega t$. The system is in the generalized multi-variable structure in (\ref{eqn:ESC_multi}) and the corresponding LBS for (\ref{eqn:ESC_BU}) is given by
\begin{equation}\label{eqn:LBS_BU}
\begin{split}
    \dot{z}_x= -\frac{k \alpha }{2}\frac{\partial f}{\partial z_x};
    \dot{z}_y= -\frac{k \alpha }{2}\frac{\partial f}{\partial z_y},
\end{split}
\end{equation}
where $\alpha = a^2$. Similar to case 1, we obtained the estimated LBS using the GEKF as described in section \ref{s: MainResults}. For the simulation, we take the following parameters $\alpha=0.5, k=2, \omega =25$ with $f(x)=x^2+y^2$. 
Figure \ref{fig:plotWithTime_BU} shows the trajectories of $x$ and $y$ coordinates of the vehicle employing ESC from \cite{BoundedUpdateKrstic} (red) and our proposed system (black). 
We verified Theorem \ref{thm:proposed_theorem} by tracking the estimated right-hand side of the LBS, i.e., $|J_{x}|$ and $|J_{y}|$ vs. time, and with a viable upper bound in the form $1/t^{1+\epsilon}$ (we choose $\epsilon =0.05$). As seen in figure \ref{fig:RHS_LBS_BU}, Theorem \ref{thm:proposed_theorem} is verified for both $x$ and $y$.

\begin{figure}[ht]
    \includegraphics[width=0.5\textwidth, keepaspectratio]{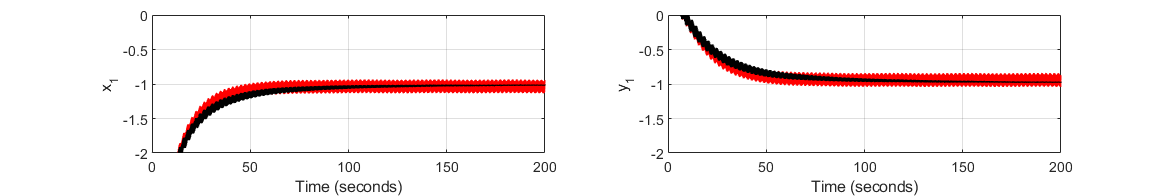}
    \caption{Trajectories of x and y coordinates of vehicle-2 employing ESC with single-integrator dynamics from \cite{DURR2013} (red) vs. our proposed system (black).}
    \label{fig: vehicle_SI}
\end{figure}

\begin{figure}[ht]
    \centering
    \includegraphics[width=0.45\textwidth, keepaspectratio]{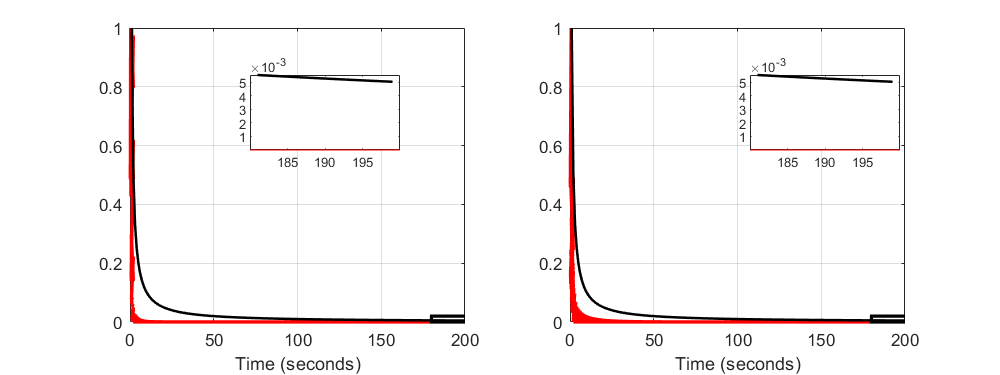}
    \caption{Verification of Theorem \ref{thm:proposed_theorem} for single-integrator dynamics by plotting $|J_{x2}|$ (red, left), $|J_{y2}|$ (red, right), and $1/t^{1+0.05}$ (black) vs. time.}
    \label{fig: vehicle2lie_SI}
\end{figure}

\noindent\textbf{Case 3.}
In this case, a problem of a three-agent vehicle system from \cite{DURR2013} is considered. We show in the appendix \ref{sec:StabilityConditionsB} that the objective function used in \cite{DURR2013} does not satisfy assumption B2, hence, this problem does not attain vanishing oscillation using the approach in \cite{VectorFieldGRUSHKOVSKAYA2018}. On the other hand, by following a similar process to the above two cases, we successfully resolved the problem with vanishing oscillations-- refer to figures \ref{fig: vehicle_SI}-\ref{fig: vehicle2lie_SI}.
%


\section{Conclusion}\label{s: Conclusion}
The proposed ESC system in this paper is based on the estimation of LBS via the novel GEKF for both stability characterization and attenuation of oscillations, and it shows promising results to control-affine ESC structures: (i) it has stability characteristics that are less complex when compared to literature (dependent on one time-dependent condition), (ii) it guarantees vanishing oscillations by design, and (iii) it resolved problems the generalized approach in \cite{VectorFieldGRUSHKOVSKAYA2018} did not succeed in solving.  

\appendix
\section{Proof of Theorem \ref{thm:proposed_theorem}}\label{sec:Proof}
Throughout the proof, $i=1,...,n$ is used to denote the $i^{th}$ component of an n-vector and we assume A1-A4 to be satisfied. It is clear that $\bm{J}$ in (\ref{eqn:LieEst}) is Riemann integrable. We fix the initial time to $t=0$ and apply $\int_0^t(\cdot) d\tau$ on both sides of (\ref{eqn:LieEst}), with $|J_i|\le {1}/{t^p} \: \forall \: t\in (t^*,\infty)$ and $p>1$, and rearranging, for all $i$, we get:

\begin{align}
    \hat{z}_i(t)-\hat{z}_{i}(0)=\int_0^t J_i d\tau =\int_0^{t^*} J_i d\tau + \int_{t^*}^t J_i d\tau.
\end{align}
Let the evaluation of the integral $\int_0^{t^*} J_i d\tau =\beta_i^*$, then,
\begin{equation}
\hat{z}_i(t)=\hat{z}_i(0)+\beta_i^*+\int_{t^*}^t J_i d\tau \label{eqn:a1}
\end{equation}
Let us consider any arbitrary two initial conditions $\bm{l}$ and $\bm{m}$ in the compact set $\mathscr{C}$ such that
\begin{align}
    \hat{z}_i(t;l_i) &=l_i+\beta^*_{i,l}   +  \int_{t_l^*}^t J_i d\tau \label{eqn:initialCond1}\\
    \hat{z}_i(t;m_i) &=m_i+ \beta^*_{i,m} +  \int_{t_m^*}^t J_i d\tau \label{eqn:initialCond2},
\end{align}
where $\beta^*_{i,l}=\int_0^{t_l^*} J_i d\tau$ and $\beta^*_{i,m}=\int_0^{t_m^*} J_i d\tau $ are some finite values. Now, from (\ref{eqn:initialCond1}) and ( \ref{eqn:initialCond2}), we get:
\begin{equation}\label{eqn:bounds_lm}
\begin{split}
      \vert l_i-m_i\vert \leq \vert \hat{z}_i(t;l_i)-\hat{z}_i(t;m_i)\vert + \vert \beta^*_{i,l}-\beta^*_{i,m}\vert +\\
      \left \vert \int_{t_l^*}^t J_i d\tau-\int_{t_m^*}^t J_i d\tau \right \vert .\\
\end{split}
\end{equation}
Let the finite quantity $\vert \beta^*_{i1}-\beta^*_{i2}\vert  \le \vert \beta^*_{i1}\vert +\vert \beta^*_{i2}\vert  \le M_1$; $M_1$ is some positive constant. With $\vert J_i\vert \le \frac{1}{t^p}$, $\forall \: t\in (t^*,\infty)$ given $p>1$,  we have
\begin{equation}
    \int_{t^*}^t J_i d\tau\leq\int_{t^*}^t \vert J_i\vert  d\tau\leq\int_{t^*}^t \frac{1}{\tau^p}d\tau. \label{eqn:time_bound}
\end{equation}
Therefore, $\left \vert \int_{t_1^*}^t J_i d\tau-\int_{t_2^*}^t J_i d\tau \right \vert $ will have some finite bound as
$\left \vert \int_{t_1^*}^t J_i d\tau-\int_{t_2^*}^t J_i d\tau \right \vert  \le
\left \vert \int_{t_1^*}^t J_i d\tau\right\vert  + \left\vert \int_{t_2^*}^t J_i d\tau \right \vert 
\le
\left \vert \int_{t_1^*}^t \frac{1}{\tau^p} d\tau\right\vert  + \left\vert \int_{t_2^*}^t \frac{1}{\tau^p} d\tau \right \vert 
\le M_2$; $M_2$ is some positive constant.
Then, (\ref{eqn:bounds_lm}) can be written as
\begin{equation}\label{eqn:lm_bounded}
    \vert l_i-m_i\vert \leq \vert \hat{z}_i(t;l_i)-\hat{z}_i(t;m_i)\vert +M,
\end{equation}
where $M=M_1+M_2$. From (\ref{eqn:lm_bounded}), it is clear that, for every $\epsilon>0$, we have $\delta=f(\epsilon)=\epsilon+M>0$ such that
\begin{equation}
    \vert l_i-m_i\vert <\delta \:\: \Rightarrow  \:\: \vert \hat{z}_i(t;l_i)-\hat{z}_i(t;m_i)\vert <\epsilon, \:\: \forall \epsilon>0.  \label{cond2}
\end{equation}
The $(\epsilon,\delta)$ stability condition in \eqref{cond2} (see section 4.3 in \cite{aastrom2008feedback} or chapter 4 in \cite{khalil2002nonlinear}) proves the stability of (\ref{eqn:LieEst}) in the Lyapunov sense. Next, we show the local asymptotic stability of the equilibrium point $\hat{z}_i^*\in \mathscr{C}$ (also is the extremum point as in A3).
Since $\lim_{t\to \infty}[\int_{t^*}^t \frac{1}{\tau^p}d\tau]$ is convergent for $p>1$, then $\lim_{t\to \infty}[\int_{t^*}^t J_i d\tau]$ is also convergent per \eqref{eqn:time_bound}. We let $\lim_{t\to \infty}[\int_{t^*}^t J_i d\tau]=c_i^*$, so with \eqref{eqn:a1}, the following condition follows:
\begin{equation}
    \lim_{t\to \infty} \hat{z}_i(t)=\hat{z}_{0i}+{\beta}_i^* +c^* = k_i^*. \label{eqn:cond1}
\end{equation}
As per (\ref{eqn:cond1}), the LBS trajectory in (\ref{eqn:LieEst}) is convergent to some constant value of $k_i^*$, which means $k_i^*$ is an equilibrium point. But from assumption A3, $\hat{z}_i^*$ is the only isolated equilibrium point. Hence, $k_i^*$ must be $\hat{z}_i^*$, and all the solutions starting from arbitrary initial conditions $l_i$ and $m_i$ within the compact set $\mathscr{C}$ converge to $\hat{z}_i^*$. Thus, $\hat{z}_i^*$ is locally asymptotically stable for the estimated LBS system (\ref{eqn:LieEst}). This completes the first part of the proof of Theorem \ref{thm:proposed_theorem}.
%
%
Now, for the proof of second part of Theorem \ref{thm:proposed_theorem}, we rewrite (\ref{eqn:law}) as: $\Dot{a_i}+ \lambda_i a_i=\lambda_i J_i$. Using the integrating factor $e^{\lambda_i t}$,
\begin{align}
    \int_0^t \frac{d}{d\tau} (a_i e^{\lambda_i \tau}) d\tau &=\int_0^t \lambda_i J_i e^{\lambda_i \tau} d\tau. \label{a12}
\end{align}
Applying $|.|$ and bounds on $J_i$ and taking  $\lim\limits_{t\to \infty}$,
\begin{equation*}
\begin{split}
    0   \le \lim\limits_{t\to \infty}|a_i(t)|  \le \cancelto{0}{\lim\limits_{t\to \infty} \left [\int_0^{t^*} \lambda_i |J_i|e^{\lambda_i \tau} d\tau \right] e^{-\lambda_i t}} \\
    +\lim\limits_{t\to \infty} \frac{\int_{t^*}^t \lambda_i \frac{1}{\tau ^p} e^{\lambda_i \tau} d \tau}{e^{\lambda_i t}}
    +\cancelto{0}{ \lim\limits_{t\to \infty}a_0 e^{-\lambda_i t}}.
\end{split}
\end{equation*}
Now, by using Taylor expansion of $e^{\lambda_i t}$ in $t$, we get
\begin{equation}\label{eqn: TaylerExp}
    0\le 
    \lim\limits_{t\to \infty} |a_i(t)|
    \le \lambda_i \lim\limits_{t\to \infty} 
    \frac{\int_{t^*}^t \frac{1}{\tau^p } (1+ \lambda_i \tau + \frac{(\lambda_i \tau) ^2}{2!} + ...) d\tau }{1+\lambda_i t+\frac{(\lambda_i t)^2}{2!}+...}.
\end{equation}
Let the value of integral in (\ref{eqn: TaylerExp}) evaluated at $t^*$ be $c$, then by writing (\ref{eqn: TaylerExp}) in a polynomial coefficient form,
\begin{equation*}
    0\le \lim\limits_{t\to \infty} |a_i(t)| \le \lambda_i \lim\limits_{t\to \infty} 
    \frac{c+a_1 t^{1-p}+ a_2t^{2-p}+...}{c_0 + c_1 t+ c_2 t^2 + ...}.
\end{equation*}
Finally, using the Squeeze theorem
\begin{align}\label{eqn:converging_a}
    \lim\limits_{t\to \infty} |a_i(t)|=0.
\end{align}
This proves the second claim of {Theorem \ref{thm:proposed_theorem}}.
%
%
%
Next, we prove the third claim of Theorem \ref{thm:proposed_theorem}. First, we note here that the proof of asymptotic stability of a class of estimated LBSs (\ref{eqn:LieEst}) includes the particular case of the exact corresponding LBS (\ref{eqn:Lie_multi}); this follows by the fact that $\eta_i(t) =0$ satisfies A4 -- leading to the exact LBS (\ref{eqn:Lie_multi}) -- being included in the proof of first claim of Theorem \ref{thm:proposed_theorem}. As a result, Theorem \ref{thm:esc_lbs} applies and by letting $\bm{x}_1(t)$ be solution of the ESC system in (\ref{eqn:generalizedSystem}) but with constant amplitude $a(t)=a_0$ as in (\ref{eqn:ESC_back}), then we have for a fixed $\omega$, for every $\epsilon >0$, $\exists$ $\delta(\epsilon) >0$ such that  $|{x}_i(0)-{x}_i^*|<\delta$ with $|{x}_{1i}(t)-{x}_i^*|<\epsilon $, where ${x}_i^*={z}_i^*$ is the equilibrium point of (\ref{eqn:generalizedSystem}) such that $\lim_{t \rightarrow \infty} \hat{z}_i(t)={z}_i^*={x}_i^*$. Now, let $\bm{x}_2(t)$ be the solution of the proposed ESC system in (\ref{eqn:generalizedSystem}) when $a(t)$ is time-varying. For both systems $\bm{x}_1(0)=\bm{x}_2(0)=\bm{x}(0)$ and $\bm{x}_1^*=\bm{x}_2^*=\bm{x}^*$. We need to show that (\ref{eqn:generalizedSystem}) is also practically asymptotically stable similar to (\ref{eqn:ESC_back}), i.e., the ($\epsilon$,$\delta$) condition satisfied above on $\bm{x}_1(t)$ is satisfied on $\bm{x}_2(t)$.
Now for fixed $\omega$, we have
\begin{equation}\label{eqn:constant_a}
    {x}_{1i}(t)={x}_{i}(0)+\int_0^t \sqrt{\omega}a_0 \hat{u}_{1i} b_{1i} d\tau + \int_0^t \sqrt{\omega}a_0 \hat{u}_{2i} b_{2i}d\tau.
\end{equation}
Equation (\ref{eqn:constant_a}) can be rewritten as 
\begin{equation*}
\begin{split}
    &\big|\int_0^t \sqrt{\omega}a_0 \hat{u}_{1i} b_{1i}d\tau + \int_0^t \sqrt{\omega}a_0 \hat{u}_{2i} b_{2i}d\tau\big| \\ &=|({x}_{1i}(t)-x_{i}^*)-({x}_{i}(0)-x_{i}^*)|\\
    &\le |x_{1i}(t)-x_i^*|+|x_i(0)-x_i^*|
     \le \delta(\epsilon) + \epsilon.
\end{split}
\end{equation*}

\noindent Then, for $\forall \epsilon>0$ we let $|\int_0^t \sqrt{\omega}a_0 \hat{u}_{1i} b_{1i}d\tau| + |\int_0^t \sqrt{\omega}a_0 \hat{u}_{2i} b_{2i}d\tau| \le \delta(\epsilon)+\epsilon= K_i(\epsilon)$. We also have:
\begin{equation}\label{eqn:variable_a}
    {x}_{2i}(t)={x}_{i}(0)+\int_0^t \sqrt{\omega}a_i(\tau) \hat{u}_{1i} b_{1i}d\tau + \int_0^t \sqrt{\omega} a_i(\tau) \hat{u}_{2i} b_{2i}d\tau,
\end{equation}
Subtracting (\ref{eqn:constant_a}) from (\ref{eqn:variable_a}), we get:
\begin{multline}\label{eqn:diff_a}
    x_{2i}(t)-x_{1i}(t)=\sqrt{\omega} \Bigl(\int_0^t a_i(\tau)\hat{u}_{1i} b_{1i}d\tau+\\
    \int_0^t a_i(\tau) \hat{u}_{2i} b_{2i}d\tau \Bigr)
    - \sqrt{\omega} \Bigl(\int_0^t  a_{0} \hat{u}_{1i} b_{1i}d\tau+ \int_0^t a_0 \hat{u}_{2i} b_{2i}d\tau \Bigr). 
\end{multline}
From (\ref{eqn: TaylerExp}) and (\ref{eqn:converging_a}), $a(t)$ is bounded for all $t \implies \exists K_0>0$ such that $|a_i(t)|\le K_0 a_0$. Now, from (\ref{eqn:diff_a})
\begin{equation}
    |x_{2i}(t)-x_{1i}(t)|\le (K_0 +1)K_i(\epsilon)
\end{equation}
Now,
\begin{align}
    x_{2i}(t)-x_i^*&=(x_{2i}(t)-x_{1i}(t))+x_{1i}(t)-x_i^*\\
    |x_{2i}(t)-x_i^*| & \le (K_0+1)K_i(\epsilon)+\epsilon =\epsilon^*
\end{align}
Thus, for ever $\epsilon^*>0$, $\exists$ $\delta>0$ such that we have $|x_i(0)-x_i^*| < \delta$ with $|x_{2i}(t)-x_i^*| <\epsilon^* $ for a given $\omega$. This completes the proof of Theorem \ref{thm:proposed_theorem}.
\section{Stability conditions in \cite{VectorFieldGRUSHKOVSKAYA2018} vs. case 3.}\label{sec:StabilityConditionsB}
 For a domain $D\subseteq \mathbb{R}^n$, refer \cite[section 3]{VectorFieldGRUSHKOVSKAYA2018}, condition B1 below is needed for practical stability whereas condition B2 is further needed for asymptotic stability (vanishing oscillations):
\begin{enumerate}[label=B\arabic*.]
    \item
    There exist constants $ \gamma_1, \gamma_2, \kappa_1, \kappa_2, \mu$ and $m_1 \ge 1$ such that for all $\bm{x} \in D$, 
    \begin{equation*}{\label{eqn:grushExample}}
        \begin{split}
    \gamma_1 ||\bm{x}-\bm{x}^*||^{2m_1} &\le \tilde{f}(\bm{x}) \le \gamma_2 ||\bm{x}-\bm{x}^*||^{2m_1},\\
    \kappa_1\Tilde{f}(\bm{x})^{2- \frac{1}{m_1}} &\le ||\nabla f(\bm{x})||^2 \le \kappa_2 \tilde{f}(\bm{x})^{2-\frac{1}{m_1}},\\
    \left | \left|\frac{\partial^2f(\bm{x})}{\partial \bm{x}^2}\right | \right| &\le \mu \tilde{f}(\bm{x})^{1-\frac{1}{m_1}},
        \end{split}
    \end{equation*}
    where $||.||$ is the $l_2$ norm, and $\tilde{f}(\bm{x})=f(\bm{x})-f^*$.
    \item
    The functions $b_{si}(\tilde{f}(\bm{x}))$ are Lipschitz on each compact set from $D$, and
    \begin{equation*}
        \begin{split}
            \alpha_1 \tilde{f}^{m_2}(x) \le b_{0i}(\tilde{f}(x))\le \alpha_2 \tilde{f}^{m_2}(x)\\
            |b_{si}(\tilde{f}(x))| \le M\tilde{f}^{m_3}(x)\\
            ||L_{b_{ql}}L_{b_{pj}}b_{si}(\tilde{f}(\cdot))|| \le H \tilde{f}^{m_4}(x)
        \end{split}
    \end{equation*}
    for all $x\in D$, $s,p,q=\overline{1,2}, i,j,l=\overline{1,n}$ with $m_2 \ge 1/m_1-1, m_3=(m_2+1)/2, m_4=3(1+m_2)/2-1/m_1$, and some $\alpha_1,\alpha_2, M>0, H\ge0$.
\end{enumerate}
Now, we check condition B2 which is needed for guaranteeing vanishing oscillations of the control input at the extremum point for the multi-agent system with single integrator dynamics in \cite[section 4]{DURR2013}, which we solved in case 3 (section \ref{s:Simulation}) successfully with vanishing oscillations. Without any loss of generality, we limit our analysis here to the objective function (map) assigned to vehicle-3 of the multi-agent system. We analyze said objective function vs. the bounds in condition B2. The expression for the objective function for vehicle-3 is $f_3 = -(x_3+1)^2/2-3(y_3-1)^2/2+10$. The extremum point is $x_3^*=-1$ and $y_3^* =1$, which corresponds to the optimal value of the cost function as $f_3^*=10$. Thus, $\tilde{f}_3=f_3-f_3^*=-(x_3+1)^2/2-3(y_3-1)^2/2$. Now, each element of the vector field $\bm{b}$ are
\begin{equation*}
    \begin{split}
    b_{11}(\tilde{f_3})&=c_3(f_3(\bm{x})-x_{3e}h_3),\: b_{21}(\tilde{f_3})=a_3,\\
    b_{12}(\tilde{f_3})&=a_3,\: b_{22}(\tilde{f_3})=-c_3(f_3(\bm{x})-x_{3e}h_3),\\
    b_{13}(\tilde{f_3})&=0,\quad b_{23}(\tilde{f_3})=0.
    \end{split}
\end{equation*}
As per condition B2, the following inequality should be satisfied for all $x\in D$
\begin{equation*}
    |b_{si}(\tilde{f}(x))| \le M\tilde{f}^{m_3}(x),
\end{equation*}
where $s\in \{1,2\}$, and $i \in \{1,2,3\}$. Let us analyze $b_{12}(\tilde{f_3})=a_3$ as an element of the vector field, and the point $x=-1,y=1$ as the point of reference. Since $M\tilde{f}^{m_3}=0$ at the reference point for any $M>0$ and $m_3$, the inequality becomes $|a_3| \le 0$. However, $a_3=constant\:>0$ is the amplitude of the input signal and is defined to be a positive constant \cite{DURR2013}, which establishes a contradiction. Thus, the condition B2 is not satisfied for vehicle-3 in the multi-agent system \cite{DURR2013}.

\bibliographystyle{IEEEtran}
\bibliography{bibliography}

\begin{thebibliography}{10}
\providecommand{\url}[1]{#1}
\csname url@samestyle\endcsname
\providecommand{\newblock}{\relax}
\providecommand{\bibinfo}[2]{#2}
\providecommand{\BIBentrySTDinterwordspacing}{\spaceskip=0pt\relax}
\providecommand{\BIBentryALTinterwordstretchfactor}{4}
\providecommand{\BIBentryALTinterwordspacing}{\spaceskip=\fontdimen2\font plus
\BIBentryALTinterwordstretchfactor\fontdimen3\font minus
  \fontdimen4\font\relax}
\providecommand{\BIBforeignlanguage}[2]{{%
\expandafter\ifx\csname l@#1\endcsname\relax
\typeout{** WARNING: IEEEtran.bst: No hyphenation pattern has been}%
\typeout{** loaded for the language `#1'. Using the pattern for}%
\typeout{** the default language instead.}%
\else
\language=\csname l@#1\endcsname
\fi
#2}}
\providecommand{\BIBdecl}{\relax}
\BIBdecl

\bibitem{KRSTICMain}
M.~Krstić and H.-H. Wang, ``Stability of extremum seeking feedback for general
  nonlinear dynamic systems,'' \emph{Automatica}, vol.~36, no.~4, pp. 595--601,
  2000.

\bibitem{Maggia2020higherOrderAvg}
M.~Maggia, S.~A. Eisa, and H.~E. Taha, ``On higher-order averaging of
  time-periodic systems: reconciliation of two averaging techniques,''
  \emph{Nonlinear Dynamics}, vol.~99, no.~1, pp. 813--836, Jan 2020.

\bibitem{ariyur2003real}
K.~B. Ariyur and M.~Krstic, \emph{Real-time optimization by extremum-seeking
  control}.\hskip 1em plus 0.5em minus 0.4em\relax John Wiley \& Sons, 2003.

\bibitem{scheinker2017modelfreebook}
A.~Scheinker and M.~Krsti{\'c}, \emph{Model-free stabilization by extremum
  seeking}.\hskip 1em plus 0.5em minus 0.4em\relax Springer, 2017.

\bibitem{ESCTracking}
C.~Zhang, A.~Siranosian, and M.~Krstic, ``Extremum seeking for moderately
  unstable systems and for autonomous vehicle target tracking without position
  measurements,'' in \emph{2006 American Control Conference}, 2006, pp. 6
  pp.--.

\bibitem{DURR2013}
H.-B. Dürr, M.~S. Stanković, C.~Ebenbauer, and K.~H. Johansson, ``Lie bracket
  approximation of extremum seeking systems,'' \emph{Automatica}, vol.~49,
  no.~6, pp. 1538--1552, 2013.

\bibitem{scheinker2016boundedDither}
A.~Scheinker and D.~Scheinker, ``Bounded extremum seeking with discontinuous
  dithers,'' \emph{Automatica}, vol.~69, pp. 250--257, 2016.

\bibitem{BoundedUpdateKrstic}
A.~Scheinker and M.~Krstić, ``Extremum seeking with bounded update rates,''
  \emph{Systems \& Control Letters}, vol.~63, pp. 25--31, 2014.

\bibitem{ExpStabSUTTNER2017}
R.~Suttner and S.~Dashkovskiy, ``Exponential stability for extremum seeking
  control systems,'' \emph{IFAC-PapersOnLine}, vol.~50, no.~1, pp.
  15\,464--15\,470, 2017, 20th IFAC World Congress.

\bibitem{eisa2023}
S.~A. Eisa and S.~Pokhrel, ``Analyzing and mimicking the optimized flight
  physics of soaring birds: A differential geometric control and extremum
  seeking system approach with real time implementation,'' \emph{SIAM Journal
  on Applied Mathematics}, 2023.

\bibitem{bullo2019geometric}
F.~Bullo and A.~D. Lewis, \emph{Geometric control of mechanical systems:
  modeling, analysis, and design for simple mechanical control systems}.\hskip
  1em plus 0.5em minus 0.4em\relax Springer, 2019, vol.~49.

\bibitem{DURR2017}
H.-B. Dürr, M.~Krstić, A.~Scheinker, and C.~Ebenbauer, ``Extremum seeking for
  dynamic maps using lie brackets and singular perturbations,''
  \emph{Automatica}, vol.~83, pp. 91--99, 2017.

\bibitem{GRUSHKOVSKAYA2017}
V.~Grushkovskaya, H.-B. Dürr, C.~Ebenbauer, and A.~Zuyev, ``Extremum seeking
  for time-varying functions using lie bracket approximations,''
  \emph{IFAC-PapersOnLine}, vol.~50, no.~1, pp. 5522--5528, 2017, 20th IFAC
  World Congress.

\bibitem{VectorFieldGRUSHKOVSKAYA2018}
V.~Grushkovskaya, A.~Zuyev, and C.~Ebenbauer, ``On a class of generating vector
  fields for the extremum seeking problem: Lie bracket approximation and
  stability properties,'' \emph{Automatica}, vol.~94, pp. 151--160, 2018.

\bibitem{taha2021vanishing}
M.~Abdelgalil and H.~Taha, ``Lie bracket approximation-based extremum seeking
  with vanishing input oscillations,'' \emph{Automatica}, vol. 133, p. 109735,
  2021.

\bibitem{AttenuatedOscillaiion2021}
D.~Bhattacharjee and K.~Subbarao, ``Extremum $ $ seeking control with
  attenuated steady-state oscillations,'' \emph{Automatica}, vol. 125, p.
  109432, 2021.

\bibitem{isidori1985nonlinear}
A.~Isidori, \emph{Nonlinear control systems: an introduction}.\hskip 1em plus
  0.5em minus 0.4em\relax Springer, 1985.

\bibitem{Chicka2006}
D.~F. Chichka, J.~L. Speyer, C.~Fanti, and C.~G. Park, ``Peak-seeking control
  for drag reduction in formation flight,'' \emph{Journal of Guidance, Control,
  and Dynamics}, vol.~29, no.~5, pp. 1221--1230, 2006.

\bibitem{ESC_KF}
G.~Gelbert, J.~P. Moeck, C.~O. Paschereit, and R.~King, ``Advanced algorithms
  for gradient estimation in one- and two-parameter extremum seeking
  controllers,'' \emph{Journal of Process Control}, vol.~22, no.~4, pp.
  700--709, 2012.

\bibitem{grushkovskaya2021extremum}
V.~Grushkovskaya and C.~Ebenbauer, ``Extremum seeking control of nonlinear
  dynamic systems using lie bracket approximations,'' \emph{International
  Journal of Adaptive Control and Signal Processing}, vol.~35, no.~7, pp.
  1233--1255, 2021.

\bibitem{aastrom2008feedback}
K.~J. {\AA}str{\"o}m and R.~M. Murray, ``Feedback systems: an introduction for
  scientists and engineers,'' 2009, Version v2.10b.

\bibitem{khalil2002nonlinear}
H.~K. Khalil, ``Nonlinear systems third edition,'' \emph{Patience Hall}, vol.
  115, 2002.

\end{thebibliography}
\end{document}